
\magnification1200
\input amssym.def
\input amssym.tex

\hsize=16truecm
\baselineskip=16truept plus.3pt minus .3pt

\font\secbf=cmb10 scaled 1200
\font\eightrm=cmr8
\font\sixrm=cmr6

\font\eighti=cmmi8

\font\sixi=cmmi6
\skewchar\eighti='177 \skewchar\sixi='177

\font\eightsy=cmsy8
\font\sixsy=cmsy6
\skewchar\eightsy='60 \skewchar\sixsy='60

\font\eightit=cmti8

\font\eightbf=cmbx8
\font\sixbf=cmbx6

\let\sc=\tensc

\font\eightsc=cmcsc10 scaled 800
\font\secbf=cmb10 scaled 1200
\font\subsecfont=cmb10 scaled \magstephalf
\font\amb=cmmib10

\font\ambi=cmmib10 scaled 700

\newfam\mbfam 

\textfont\mbfam\amb \scriptfont\mbfam\ambi


\def\aa{\def\rm{\fam0\eightrm}%
  \textfont0=\eightrm \scriptfont0=\sixrm \scriptscriptfont0=\fiverm
  \textfont1=\eighti \scriptfont1=\sixi \scriptscriptfont1=\fivei
  \textfont2=\eightsy \scriptfont2=\sixsy \scriptscriptfont2=\fivesy
  \textfont3=\tenex \scriptfont3=\tenex \scriptscriptfont3=\tenex
  \def\sc{\eightsc}
  \def\it{\fam\itfam\eightit}%
  \textfont\itfam=\eightit
  \def\bf{\fam\bffam\eightbf}%
  \textfont\bffam=\eightbf \scriptfont\bffam=\sixbf
   \scriptscriptfont\bffam=\fivebf
  \normalbaselineskip=9.7pt
  \setbox\strutbox=\hbox{\vrule height7pt depth2.6pt width0pt}%
  \normalbaselines\rm}

\def\Proof{\vskip12pt\noindent{\bf Proof.} }

\def\Proposition#1{\vskip12pt\noindent{\bf Proposition #1}}

\def\m@th{\mathsurround=0pt}

\def\cc#1{\hbox to .89\hsize{$\displaystyle\hfil{#1}\hfil$}\cr}
\def\lc#1{\hbox to .89\hsize{$\displaystyle{#1}\hfill$}\cr}
\def\rc#1{\hbox to .89\hsize{$\displaystyle\hfill{#1}$}\cr}

\def\eqal#1{\null\,\vcenter{\openup\jot\m@th
  \ialign{\strut\hfil$\displaystyle{##}$&&$\displaystyle{{}##}$\hfil
      \crcr#1\crcr}}\,}

\def\section#1{\vskip 22pt plus6pt minus2pt\penalty-400
        {{\secbf
        \noindent#1\rightskip=0pt plus 1fill\par}}
        \par\vskip 12pt plus5pt minus 2pt
        \penalty 1000}

\def\subsection#1{\vskip 20pt plus6pt minus2pt\penalty-400
        {{\subsecfont
        \noindent#1\rightskip=0pt plus 1fill\par}}
        \par\vskip 8pt plus5pt minus 2pt
        \penalty 1000}

\def\subsubsection#1{\vskip 18pt plus6pt minus2pt\penalty-400
        {{\subsecfont
        \noindent#1}}
        \par\vskip 7pt plus5pt minus 2pt
        \penalty 1000}

\def\centerlast#1{\begingroup \leftskip=0pt plus 1fil\rightskip=0pt plus -1fil
\parfillskip=0pt plus 2fil
\parindent 0pt
\par #1\par\endgroup}
\def\\{\hfill\break}

\def\kwadrat{\null\ \hfill\null\ \hfill$\square$}
\def\mida#1{{{\null\kern-4.2pt\left\bracevert\vbox to 6pt{}\!\hbox{$#1$}\!\right\bracevert\!\!}}}
\def\midy#1{{{\null\kern-4.2pt\left\bracevert\!\!\hbox{$\scriptstyle{#1}$}\!\!\right\bracevert\!\!}}}

\def\divv{{{\rm div}\,}}

\def\today{${\scriptscriptstyle\number\day-\number\month-\number\year}$}
\footline={{\hfil\rm\the\pageno\hfil${\scriptscriptstyle\rm\jobname}$\ \ \today}}

\def\ifnextchar#1#2#3{\bgroup
  \def\reserveda{\ifx\reservedc #1 \aftergroup\firstoftwo
    \else \aftergroup\secondoftwo\fi\egroup{#2}{#3}}%
  \futurelet\reservedc\ifnch
  } 
\def\ifnch{\ifx \reservedc \sptoken \expandafter\xifnch
      \else \expandafter\reserveda
      \fi} 
\def\firstoftwo#1#2{#1}
\def\secondoftwo#1#2{#2} 
\def\tempswafalse{\let\iftempswa\iffalse}
\def\tempswatrue{\let\iftempswa\iftrue} 

\def\cite{\ifnextchar [{\tempswatrue\citea}{\tempswafalse\citeb}}
\def\citea[#1]#2{[#2, #1]}
\def\citeb#1{[#1]}\phantom{]}

\def\\{\hfil\break}
\def\C{{\Bbb C}}

\def\N{{\Bbb N}}
\def\R{{\Bbb R}}
\def\Z{{\Bbb Z}}

\def\dH{{\dot H}}
\def\dX{{\dot X}}
\def\dZ{{\Z^d \setminus \{ 0 \}}}
\def\divv{{\rm div}}
\def\span{{\rm span}}
\def\esssup{{\rm ess \, sup}}
\def\nad#1#2{{\mathop{#1}\limits^{#2}\null\!}}
\def\longlongrightarrow{{-\!\!\!-\!\!\!-\!\!\!\longrightarrow}}

\begingroup
\catcode`\&=\active
\gdef\url{%
  \begingroup%
    \def\urlaux##1{##1\endgroup}%
    \catcode`\&=\active
    % Catcode of & that is important here is when \url is defined:
    \def&{\penalty0\char`\&}% 
    \urlaux}
\endgroup

\centerlast{\secbf Quantitative robustness of regularity for 3D Navier-Stokes system in 
$\dH^\alpha$-spaces}

\vskip1cm
\centerline{\bf Jan Burczak$^1$ and W. M. Zaj\c aczkowski$^{1,2}$}
\vskip1cm
\item{$^1$} Institute of Mathematics, Polish Academy of Sciences,\\ 
\'Sniadeckich 8, 00-950 Warsaw, Poland,\\
e-mail: jb@impan.pl;\\
\item{$^2$} Institute of Mathematics and Cryptology, Cybernetics Faculty,\\ 
Military University of Technology, Kaliskiego 2, 00-908 Warsaw,\\ 
Poland\\

\vskip1cm

\noindent
{\bf Abstract.} 
We present stability and regularity results for the $3$D incompressible Navier-Stokes system in 
a periodic box, in $\dH^\alpha$ spaces, with $\alpha\in\big[{1/2},1\big]$. A special attention is paid to obtaining quantitative results,   i.e.  ones with explicit or at least computable constants, and to scaling.
\vfil\eject

\section{1. Introduction}
We study stability of a regular solution $u$ of a $3$D Navier-Stokes system in the periodic cube $Q_L := [0,L]^3$. Namely, let us
fix an $L$-periodic pair $(u,p)$ solving\footnote{$^{(1)}$}{For rigorous 
definition of the solution and for presentation of underlying function space, please refer to Section 2.} in $Q_L\times[0,T]$
$$\left\{\eqal{
&u,_{t}+u\cdot\nabla u-\nu\Delta u+\nabla p=f,\cr
&\divv u=0,\cr
&u(0)=u_0,\cr}\hskip2cm (NS_{f,u_0})\right.
$$
where $\nu$ is a positive  parameter (viscosity) and  $f$ is a given function (external force). We assume that 
$u$ is more regular than a weak solution;
more precisely, that it is additionally an $\alpha$-strong solution,   i.e. 
$u\in L^\infty(\dH^\alpha)\cap L^2(\dH^{1+\alpha})$ for 
$\alpha\in\big[{1/2},1\big]$. Next, let us consider a weak solution $v$ to $(NS_{g,v_0})$ and ask what are the 
conditions on differences of data of $u$ and $v$,   i.e.  on $|f-g|$ and $|u_0-v_0|$, that allow $v$ to inherit $\alpha$-strong regularity of $u$.

This problem is referred to as a problem of stability of strong solutions or, in a more debonair manner, as a problem of robustness of regularity. It can be seen as an intermediary step between the easily accessible small-data-regularity of solutions to the $3$D Navier-Stokes system and their unknown regularity {\it in the large}, which is one of the Millennium Problem of the Clay Mathematics Institute.  Namely, one may hope that it is possible, firstly, to cover an entire space of initial data with a net of initial data that give rise to regular solutions and, next, to conclude the regularity for {\it every} initial  datum by a stability result around points of this net. Some of such programs aimed at obtaining regularity by stability are so-called  {\it schemes for numerical verification of regularity}, see Mar\'\i n-Rubio, Robinson and Sadowski \cite{M-RRS} and its references (for more on this, compare also the concluding section of this note).

In order to make a stability result useful in practice for further attempts to attack the regularity problem, it needs to contain clear dependences on parameters and constants (which especially important for numerics). This was our initial motivation to prove the main result of this note, namely Theorem 1, subsection 1.2. It refines and generalizes Theorem 1 of \cite{M-RRS}. Our Theorem 1 provides explicit, quantitative dependencies on parameters and constants. Moreover, it takes into account  scaling-related phenomena. For the next planned steps of our studies of regularity-via-stability, that shall originate in this note, please refer to its concluding section. There, we suggest also a new approach to devise a scheme aimed at {\it numerical falsification of regularity conjecture of solutions to the $3$D Navier-Stokes system}.

As a byproduct (or, more precisely, as a needed ingredient to prove  Theorem 1) we obtain also a global-in-time regularity result for 
small data and a regularity result for a small existence time (respectively, Theorems 2 and 3 in subsection 1.2). In the former, in addition to the standard blowup characterization of the maximal time of existence, we provide also a caloric characterization. These regularity results are also quantitative, in the sense of explicitly computed constants and smallness parameters.

\subsection{1.1. Current state of research in stability of strong solutions to the Navier-Stokes 
system}

A numerous variants of the stability problem have been a subject of  intensive research. The following, extremely brief presentation merely recalls  the most common approaches. In particular, we do not dwell into a diversity of the considered domains or boundary conditions.

One often proves stability separately for each special regular solution separately. Namely, one fixes a concrete `special'  regular solution (for instance a two-dimensional one, an axially symmetric without swirl one, helicoidal one etc.) and shows that any perturbed solution (i.e.  one with data close to the fixed, regular one) stays regular. A more general approach consists in taking any solution from a given class of regular solutions and showing that its perturbations remain regular. 
Compare for instance Iftimie \cite{I}, Mucha \cite{M1} (for two dimensional special solutions), Zaj\c{a}czkowski \cite{Za} (where the stability problem around a given linear combination of two-dimensional and axially symmetric solutions is considered), Bardos et als. \cite{BLNNT} (for all three mentioned types),  Zaj\c{a}czkowski and Zadrzy\'nska \cite{ZaZa} as well as their references.

Finally, one can simply consider an arbitrary regular solution (without indicating its construction or class of special solutions that it belongs to) and show that its perturbations are regular. In this context compare for instance: DaVeiga and Secchi \cite{dVS} for $L^p$ spaces approach; Auscher, Dubois and Tchamitchian \cite{ADT} for $BMO^{-1}$ spaces, see also Koch and Tataru \cite{KT} for the latter. 

Let us finally recall Mucha \cite{M2}, the monograph by Chemin et als. \cite{CDGG} and the classical one by Constantin and Foias \cite{CF}, because our considerations are close to them: to \cite{M2} in context of providing a stability result in a periodic setting, to \cite{CDGG} in context of similar approach to regularity and used function spaces, and to \cite{CF} in context of a functional setting and a special attention to scaling.

\subsection{1.2. The results}
Recall that $\nu>0$ is a viscosity parameter, $Q_L := [0,L]^3$ is our periodic domain and that $\alpha$-strong solution to a Navier-Stokes system is such  weak solution, that additionally belongs to
$L^\infty(\dH^\alpha)\cap L^2(\dH^{1+\alpha})$. We denote the Fourier-series-based norm in $\dH^\alpha (Q_L)$ with $| \cdot |_{\alpha, L}$. All the needed (standard) definitions has been shifted to Section 2.
Let us fix 
$$\eqal{
&K_2:= \sqrt{2} C_S(1-\alpha) C_S(1) C_S\bigg(\alpha-{1\over2}\bigg) \bigg({2\pi \over L}\bigg)^{-1},\cr
&K_3:=\varepsilon_1^{-3}{27\over128}  (2 \pi)^{ { -12} } C^4_S(1-\alpha) C^4_S(1)    \bigg({2\pi\over L}\bigg)^{2 (1-2 \alpha)}  \bigg[1
+  C_S\bigg(\alpha-{1\over2}\bigg) (2 \pi)^{ { \alpha -2}}  \bigg]^4,\cr
&K_4:={1\over4 \varepsilon_2} \bigg({2\pi\over L}\bigg)^{-2},}
$$
where  $C_S (\beta)$ denotes a constant of the $2 \pi$-normalized, $3$D Sobolev-Poincar\'e inequality,   i.e. 
$$
|f|_{L^{\beta^*} (Q_{2 \pi}) } \le C_S (\beta) |f|_{\beta, {2\pi}},
$$
where $\beta \in [0,2), \; \beta^* := {6 \over 3 - 2\beta}$. Since null-mean-value functions are involved, the lower order terms for the r.h.s. of the inequality above are superfluous.

Observe that for $K_3$ becomes invariant upon $L$-scaling for $\alpha = {1 \over 2}$.

\subsubsection{1.2.1 The stability result}

\proclaim Theorem 1. (Robustness of regularity).
Let us choose $T_*> 0$, $\alpha \in [{1 /2}; 1]$ and data  
$$u_0, v_0 \in \dH_\divv^\alpha(Q_L), \quad f, g\in L^2(0,T_*;\dH_\divv^{\alpha-1}(Q_L)). $$
Assume that $u$ is an $\alpha$-strong solution to $(NS_{f,u_0})$ 
with its time of existence $T_*$. Given
any positive $\bar\nu, \varepsilon_1, \varepsilon_2$ such that 
$$\bar\nu + \varepsilon_1+\varepsilon_2<\nu,$$
\noindent
 every Leray-Hopf weak solution $v$ that starts close to $u$ and that 
has a similar forcing is also an $\alpha$-strong solution.  \\
\noindent
More precisely, let us fix any $T<T_*$. Under the proximity assumption for the data
$$
\bigg(|u_0-v_0|_{\alpha, L}^2+ K_4 \int_0^T 
|f-g|^2_{\alpha-1, L} (t) \, dt\bigg) \; {\rm e}^{ K_3 \int_0^T
|\nabla u (t)|_{L^{3\over2-\alpha} (Q_L)}^4 \, dt} 
<\bigg({\bar\nu\over K_2}\bigg)^2
\eqno(A1)
$$
$v$ is an $\alpha$-strong solution to 
$(NS_{g,v_0})$ with its time of existence $T_*(g, v_0)>T$. Moreover, $v$ is close to $u$ according to the following formula
$$
\sup_{t\in[0,T]}|u-v|_{\alpha, L}^2 (t) + (\nu - (\bar\nu + \varepsilon_1+\varepsilon_2)) 
\int_0^T|u-v|_{\alpha+1, L}^2 (t) \,dt\le  \bigg({\bar\nu\over K_2}\bigg)^2.
\eqno(P1)
$$

\noindent
The proof of Theorem 1 is provided in Section 3.

 Observe that  $\intop_0^T|\nabla u(t)|_{L^{3\over2-\alpha}}^4 dt < \infty$, 
needed in the proximity assumption $(A1)$, is given automatically by the fact that $u$ is the $\alpha$-strong 
solution, since by interpolation
$$
L^\infty(\dH^\alpha)\cap L^2(\dH^{1+\alpha})\hookrightarrow
L^4(W^{1,{3\over{2-\alpha}}}).
$$
To clarify this point quantitatively, let us state
\proclaim Corollary 1. For validity of the proximity assumption (A1) it suffices
$$\eqal{
\bigg(|u_0-v_0|_{\alpha, L}^2+ K_4 \int_0^T 
|f-g|^2_{\alpha-1, L} (t) \, dt\bigg) \times \cr 
{\rm e}^{ K_3 C^4_I (\alpha, T) \big({4\pi^2\over L^2}\big)^{2(\alpha-1 )} |u|^{2}_{L^\infty(0,T; |\cdot|_{\alpha, {L}})} |u|^{2}_{L^2(0,T;  |\cdot|_{1+ \alpha, L})}
} 
<\bigg({\bar\nu\over K_2}\bigg)^2,}
\eqno(A2)
$$
where $C_I (\alpha, T) $  comes from Definition 2 in subsection 2.5.

\noindent
The proof of Corollary 1 is given in Section 3.

 Let us remark that having a qualitative $\dH^{1\over2}$ stability result  (see for instance Theorem 3.5 in [CDGG]), one  can immediately state a stability result in $\dH^{\alpha}, \, \alpha > 1/2$ via a Ladyzhenskaya-Prodi-Serrin-type condition, but it would be a merely qualitative one (without any control of involved constants and parameters), hence useless for practical applications.

\subsection{1.2.2 The regularity results}
\subsubsection{Local-in-time regularity and uniqueness}
In order to show Theorem 1, we need the following theorem on local-in-time existence of strong solutions  and on their uniqueness

\proclaim Theorem 2. (local-in-time $\alpha$-strong solutions). 
Given $T>0$, $\alpha \in [{1/2}, 1]$ and data
 $$u_0\in \dH_\divv^\alpha(Q_L), \;
f\in L^2(0,T;\dH_\divv^{\alpha-1}(Q_L)),$$
 there is $ T_*(u_0,f) \in (0, T]$ such that there 
exists a 
$$C([0,  T_*(u_0,f));
\dH_\divv^\alpha(Q_L))\cap L^2(0,T_*(u_0,f);\dH_\divv^{1+\alpha}(Q_L))$$
solution to $(NS_{f, u_0})$. Moreover, such solution is unique among Leray-Hopf weak solutions to $(NS_{f, u_0})$. \\
\noindent
Finite $T_*(u_0,f)$ can be characterized by the blowup,   i.e. 
$$
\max_{t \le \tau} |u(t)|_{\alpha,L}
\nad\longlongrightarrow{\tau \to T_*(u_0,f)}\infty.
$$

\noindent
Observe that in Theorem 2 one obtains continuity in time of the $\dH_\divv^\alpha$ norm and not only its boundedness in time, that is in the definition of an $\alpha$-strong solution. The proof of Theorem 2 is standard. For clarity we present it in Subsection 4.2. There, we obtain also the following {\it caloric lower bound} for the time $T_*(u_0,f) $. In order to formulate it, for  $k_0 \in \dZ $ let us denote by 
$$P_{k_0}:\dH_\divv^0(Q)\to\span\{e^{2\pi ik\cdot{x\over L}}:k\le k_0\}$$
the projection on the low-frequency space (with the inequality understood componentwise)

\noindent
In the following lemma, by $u^{Lo}$ we understand the solution to the homogeneous heat system that emanates from $P_{k_0} u_0$,   i.e.  to
$$\left\{ \eqal{
&u,_{t}^{Lo} - \, \nu \Delta u^{Lo}=0,\cr
&u^{Lo}(0)=u_0^{Lo}.} \right. 
$$
Recall that $K_2$ and $K_3$ come from Theorem 1.
\proclaim Lemma 1. (Caloric lower bound for $T_*$). Let us fix
any positive $ \varepsilon_1, \varepsilon_2$ and  $k_0 \in  \dZ$ so large  that 
$$\mu := \nu - \bigg( \varepsilon_1+\varepsilon_2 + K_2  \bigg({1 \over \sqrt{2}} | u_0 - P_{k_0} u_0|_{\alpha,L} \bigg)  \bigg) > 0.$$
as well as  an arbitrary $\delta>0$ and $\sigma \in (0,1)$.
Any time ${T_0}$ that yields
$$\eqal{
& \bigg( {K_2 \over {\nu - \varepsilon_1 - \varepsilon_2-{ \sigma \mu} }} \bigg)^2 {1 \over 2 \delta}  \bigg({4\pi^2\over L^2}\bigg)^{2\alpha+1} \int_0^{T_0}   | 
u^{Lo} (\tau) \otimes u^{Lo} (\tau)|^2_{1+ \alpha, L}  \, d\tau  \cr 
&\le e^{ -K_3 \big(\int_0^{T_0}   |\nabla u^{Lo} (\tau)|_{L^{3\over2-\alpha}(Q_L)}^4  d \tau + \delta T_0 \big)} - \bigg({{(1-\sigma)\mu} \over{\nu - \varepsilon_1 - \varepsilon_2-{\sigma \mu}} }-1\bigg)^2}
 \eqno(A3)
$$
is a lower bound for  $T_*(f, u_0)$, i.e. $T_*(f, u_0)\ge T_0$. Moreover
$$
\sup_{t \le {T_0}}  {1\over2} |u (t) |_{\alpha,L}^2 + {\sigma} {4\pi^2\over L^2} \int_0^{T_0} 
|u (t) |_{\alpha+1,L}^2 \, d t \le \bigg({\nu - \varepsilon_1 - \varepsilon_2-{\sigma \mu} \over K_2}\bigg)^2.
$$

\noindent
The notion {`{caloric}'} 
indicates that our lower bound is related to the homogeneous heat system, governing $u^{Lo}$. Observe that in Lemma 1 the largest $T_0$ is related to equality on (A3). We decided to keep inequality in formulation of Lemma 1, since it is easier computable than an equality.

\subsubsection{Global-in-time regular, small solutions}
Finally, we obtain also the following  global-in-time $\alpha$-strong solutions for small data.

\proclaim Theorem 3. 
Let us fix any $T\in(0,\infty]$ and any positive $\bar\nu, \varepsilon_2$ such that 
$$\bar\nu + \varepsilon_2<\nu.$$ Assume that data $f$,  $u_0$ satisfy the following 
smallness condition
$$
\bigg(|u_0|_{\alpha, L}^2+ K_4 \int_0^T 
|f (t) |^2_{\alpha-1, L} \, dt\bigg) <\bigg({\bar\nu\over K_2}\bigg)^2.  \eqno(A4)
$$
Then $(NS_{f,u_0})$ has the $\alpha$-regular solution on $[0, T]$ with the estimate
$$
\sup_{t\in[0,T]}|u (t)|_{\alpha, L}^2+ (\nu - \bar\nu -\varepsilon_2)
\int_0^T|u|_{\alpha+1, L}^2\le  \bigg({\bar\nu\over K_2}\bigg)^2.
$$

\noindent
For the proof, see Subsection 4.3. 

Since considering certain classes of 'highly--oscilliating' initial data, giving rise to regular solutions, has gained recently serious attention, let us state
\proclaim Corollary 2. Assume that $f=0$ and $\alpha \in [{1/2}, 1]$. Given any positive: $\bar\nu$ such that 
$\bar\nu <\nu$ and $L_0$ such that ${L  \over L_0} \in \N$, if the initial datum $u_0$ is $L_0$-periodic and $$ L_0 \, |(u_0)_{|L_0}|_{\alpha, L_0} <   {2 \pi \bar\nu\over  \sqrt{2} C_S(1-\alpha) C_S(1) C_S\big(\alpha-{1\over2}\big)},$$ then $(NS_{0,u_0})$ has the $\alpha$-regular solution on $[0, \infty)$ with the estimate
$$
L \sup_{t\in[0,\infty)}|u (t)|_{\alpha, L} \le {2 \pi \bar\nu\over  \sqrt{2} C_S(1-\alpha) C_S(1) C_S\big(\alpha-{1\over2}\big)}.
$$

The proof of Corollary 2 boils down to constructing our $L$-periodic solution $u$ as a juxtaposition of  ${L  \over L_0}$  $L_0$-periodic solutions staring from $(u_0)_{|L_0}$ and its shifts.

 \section{2. Preliminaries}
Here we present the detailed setting for our problem. It is standard and 
based on \cite{CF}, Chapter 4 and \cite{Tem}. 

\subsection{2.1. Function spaces}
{\bf Homogeneous Sobolev spaces $\dH^s$.}
Let us introduce the Fourier basis 
$$\omega_{L,k}^j=e^{{2\pi i\over L}k\cdot x}b_j,\quad 
\omega_{L,k}=(\omega_{L,k}^1,\cdots,\omega_{L,k}^N),
$$
with $k\in\Z^d$, 
$x\in Q_L=[0,L]^d$ and $b_j$ is the $j$-th canonical vector of $\R^N$. 
The space
$$
\dH^s(Q_L):=\bigg\{u= \!\!\! \!\!\!    \sum_{k\in\dZ}u_k\cdot\omega_{L,k}\; | \;u_k\in\C^\N,
u_k=\bar u_{-k},
 u_0=0,\sum_{k\in\dZ} \!\!\!  |k|^{2s}|u_k|^2<+\infty \bigg\},
$$
where $s \in \R$, becomes the Hilbert space, equipped with the product
$$
\langle u, w \rangle_{s,L}:= \sum_{k\in\dZ}|k|^{2s} u_k \bar w_k
$$
that generates the norm $| \cdot |_{s, L}$. We will also use the generalized scalar product (duality formula)
$$
\langle u, w \rangle_{\alpha, \beta; L}:= \sum_{k\in\dZ}( |k|^{\alpha} u_k) ( |k|^{\beta} \bar w_k)
$$
for $\alpha, \beta \in \R$. For  $\alpha + \beta = 2s$ one has $\langle u, w \rangle_{\alpha, \beta, L} = \langle u, w \rangle_{s, L}$.

We call $\dH^s(Q_L)$ a real $(u_k^l=\bar u_{-k}^l)$, zero average ($u_0^l$ plays no role, because we sum over $\dZ$), 
fractional $(s\in\R)$, homogenous Sobolev space of periodic functions 
$(u^l(x+Le_j)=u^l(x)$ thanks to the Fourier-series-based definition). 
The homogeneity of $\dH^s(Q_L)$ follows from absence of lower-order terms in its norm.

Observe that one has the following scaling-invariance. Let us define for $u:Q_L\to\R^N$ its 
dilation $u_\delta: \, Q_{\delta L}\to\R^N$ by $u_\delta(x)=u(\delta x)$. 
The Fourier coefficients of $u$ and $u_\delta$ are identical, because 
$\omega_{L,k}(x)=\omega_{\delta L,k}(\delta x)$.

One of advantages of working with homogeneous Sobolev spaces is that for 
any $s\in\R$
$$
(\dH^s(Q_L))^*\simeq\dH^{-s}(Q_L),
$$
see \cite{Tem}. 

From now on, we work with domain and target dimensions equal $3$, 
  i.e.  $d=N=3$.

We will use also the zero-divergence subspace of $\dH^s(Q)$, i.e.
$$
\dH_\divv^s(Q):=\{u\in\dH^s(Q)|\ k\cdot u_k=0, \;\; {k\in\dZ}\}
$$
closed under the norm $|\cdot|_s$.

{\bf Lebesgue spaces $L^p$ and Sobolev spaces $W^{1, p}$.}  Additionally, we will use $L^p (\Omega)$ spaces with the integro-differential norm
$|f|_{L^p (\Omega)} := \left( \int_\Omega |f|^p \right)^{1 \over p}$ as well as $W^{1, p} (\Omega)$ spaces with the integro-differential norm
$|f|_{W^{1, p}(\Omega)} := \left( \int_\Omega (|f|^p + | \nabla f|^p) \right)^{1 \over p}$. In the last formula $\nabla f : \Omega \to \R^{3 \times 3}$ denotes the weak derivative of $f$.

\noindent
At few occasions we will need the following identity
$$
|\nabla u|_{L^2 (Q_L)} =  L^{d \over 2} \bigg({4\pi^2\over L^2}\bigg)^{1 \over 2} | u |_{1, L},
\eqno(1.1)$$
valid by a computation.

{\bf Bochner spaces.}
For a bounded (time)  interval $I$ and a Banach space $V$, space $L^p (I; V)$ consists of all strongly measurable functions $u: I \to V$ with finite 
$$
|u|^p_{L^p (I; \,V)} :=  \int_I |u (t) |_V^p dt,  \qquad p \in [1, \infty);  \qquad \quad  |u|_{L^\infty (I; \, V)} := \esssup_{ t \in I}  |u (t) |_V.
$$
By a strongly measurable function we understand $u\!:I \to V$ that can be almost everywhere approximated by step functions $u_n \!: I \to V$ with respect to the norm of the space $V$, i.e. for a.e. $t \in I$
$$
\lim_{n \to \infty} |u_n (t) - u (t) |_V \;=\; 0
$$
Analogously, for a compact interval $I$, space $C (I; V)$ consists of all  continuous functions $u\!: \overline{I} \to V$ with finite
$$
 |u|_{C (I ; V )} := \max_{ t \in  I } |u (t) |.
$$
For $I = \bigcup_{n \in \N} I_n$, $I_n$'s being compact intervals,  space $C (I; V)$ consists of all  continuous functions $u\!: I \to V$ with finite
$$
 |u|_{C (I_n; V )}
$$
for any $n$. 
 
 For some more details on Bochner spaces, one may refer to Zeidler \cite{Z}, Chapter 23 and Chapter 2 of Pokorn\'y \cite{Pok}.

\subsection{2.2. Stokes operator}

The (stationary) Stokes problem in a periodic cube $Q_L$, i.e. the problem of 
finding  for a certain $f\in H^{-1}(Q_L)$ a pair $(u,p)\in\dH_\divv^1(Q_L)\times L^2(Q_L)$ that satisfies
$$
-\Delta u+\nabla p=f,\quad \divv u=0
$$
 in $Q_L$ admits in our periodic setting the following explicit solution
$$u=\sum_{k\in\dZ}u_k \cdot \omega_{L,k}, \qquad p  =\sum_{k\in\dZ} p_k \, e^{{2\pi i\over L}k\cdot x},$$
where
$$
u^j_k=-{L^2\over4\pi^2|k|^2}\bigg(f_k-{(k\cdot f_k)k\over|k|^2}\bigg), \; {\rm for } \; j=1,2,3, \qquad
p_k={Lk\cdot f_k\over2i\pi|k|^2}.
$$

Under assumption of divergent-free forcing the pressure vanishes and the `solution mapping' $f\mapsto u$ from $\dH_\divv^0(Q_L)$ to 
$\dH_\divv^2(Q_L)$ is bijective. Hence its inverse is meaningful. 
We denote it by $A:\dH_\divv^2(Q_L)\mapsto\dH_\divv^0(Q_L)$ and refer to as 
the Stokes operator. In our case it degenerates to $-\Delta$. On the side of Fourier coefficients, $A$ is the 
multiplication with $-{4\pi^2|k|^2\over L^2}$,   i.e. 
$$
(Au)_k=-u_k{4\pi^2|k|^2\over L^2}
$$
Hence we have
$$
Au=-\sum_{k\in\dZ}\lambda_ku_k\cdot\omega_{L, k}
$$
with $\lambda_k:={4\pi^2|k|^2\over L^2}$. Consequently, via the definition of norm $|\cdot |_{s,L}$,
$$
|Au|_{0,L}={4\pi^2\over L^2}|u|_{2,L}. \eqno(2.1)
$$
Formula (2.1)  admits a generalisation that defines powers of the operator 
$A$. Namely, $A^\alpha:\dH_\divv^{2\alpha}(Q_L)\mapsto\dH_\divv^0(Q_L)$ is given as
$$
A^\alpha u :=\sum_{k\in\dZ}\lambda_k^\alpha u_k\cdot\omega_{L, k}.
$$
The formula (2.1) is thus generalised to
$$
|A^\alpha u|_{0,L}=\bigg({4\pi^2\over L^2}\bigg)^\alpha
|u|_{2\alpha,L}\quad \alpha\in\R.
\eqno(2.2)
$$
and further to
$$
|A^\alpha u|_{\beta,L}=\bigg({4\pi^2\over L^2}\bigg)^{\alpha-\gamma}
|A^\gamma u|_{\delta,L} 
\eqno(2.3)
$$
for
$$
2\alpha+\beta=2\gamma+\delta,\quad \alpha,\beta,\gamma,\delta\in\R.
$$
In view of the definition of the operator $A$, we see that 
$A^{\alpha+\beta}=A^\alpha\circ A^\beta$. 

\subsection{2.4. Weak solution to Navier-Stokes system}
\noindent
In this subsection we drop a precise control over constants, because it is superfluous here. 

\noindent
Let us denote by $B(a,b)=a\cdot\nabla b$ and choose $$f\in L^2(0,T;\dH_\divv^{-1} (Q_L)), \; u_0\in\dH_\divv^0 (Q_L).$$ The first energy inequality motivates that $u$ solving 
$(NS_{f,u_0})$ belongs to $$ L^\infty(0,T; \dH_\divv^0 (Q_L))\cap L^2(0,T;\dH_\divv^1 (Q_L)).$$ In particular $u(t)\in\dH_\divv^1 (Q_L)$ for  a.a. $t\in[0,T]$. Consequently one has $$Au=-\Delta u \in L^2(0,T;\dH_\divv^{-1} (Q_L)), \quad{\rm and} \quad
B(u,u)\in L^{4 \over 3}(0,T;\dH_\divv^{-1} (Q_L)), \eqno(2.4)$$ whereas the latter follows, via duality, from the estimate
$$
| \int_0^T \langle u(t)\cdot\nabla u(t), \varphi (t) \rangle_{0, L} \; dt | \le C|u|_{L^4 (0, T; L_3 (Q_L))} |u|_{L^2 (0, T; \dH^1 (Q_L))}
|\varphi|_{L^4 (0, T; L_6 (Q_L))}
\eqno(2.5)
$$
for a sufficiently regular, divergence-free $\varphi$. The term $|u|_{L^4 (0, T; L_3 (Q_L))}$ in (2.5)  is finite thanks to parabolic embedding following from the regularity of the first energy inequality.
Hence for $f \in L^2 (0, T;  \dH_\divv^{-1} (Q_L))$, after testing formally
$(NS_{f,u_0})$ with divergence--free, sufficiently smooth, $Q_L$-periodic $\varphi$, we obtain
$$\eqal{
&\int_0^T \langle u,_{t}\!(t),\varphi (t) \rangle_{-1,1;L} dt
=  \int_0^T  \langle f(t)-B(u(t),u(t))-\nu A(u(t)),
\varphi (t)\rangle_{-1,1;L} dt \cr
& \le C_{u,f} |\varphi|_{L^4 (0, T; \dH_\divv^{1} (Q_L))},}
\eqno(2.6)
$$
where $C_{u,f}$ denotes a finite quantity related to both (2.4) and $L^2 (0, T;  \dH_\divv^{-1} (Q_L))$ norm of $f$.
Consequently (2.6) indicates, via duality, that $u,_{t}\! \in L^{4 \over 3} (0,T; \dH_\divv^{-1}(Q_L))$. This motivates the following definition 

\proclaim Definition 1. (weak solution)
Let $$f\in L^2(0,T;\dH_\divv^{-1} (Q_L)), \quad u_0\in\dH_\divv^0 (Q_L).$$ We call 
$$u\in L^\infty(0,T; \dH_\divv^0 (Q_L))\cap L^2(0,T;\dH_\divv^1 (Q_L))$$ with the 
distributional time derivative $u,_{t}\in L^{4 \over 3}(0,T; \dH_\divv^{-1} (Q_L))$ the {
(variational) weak solution} to $(NS_{f,u_0})$ iff the formula (2.6) holds for every test function 
$\varphi\in\dH_\divv^1 (Q_L)$ at almost every $t\in[0,T]$. The initial 
condition is attained in the $C_\omega([0,T];\dH_\divv^0 (Q_L))$ sense, namely 
$$
\langle u(t), \xi \rangle_{0, L} \mathop{\longrightarrow}^{t\to0}\langle u_0, \xi \rangle_{0, L}  \quad
{\rm for\ any}\quad \xi\in\dH_\divv^0 (Q_L).
$$

\noindent
Let us  motivate the way in which the initial condition is satisfied.
Since $$u,_{t}\in L^{4 \over 3}(0,T; \dH_\divv^{-1} (Q_L)) \subset L^1(0,T;\dH_\divv^{-1}(Q_L)), $$ $u$ has a representative in $C([0,T],\dH_\divv^{-1} (Q_L))$. This information 
together with $u\in L^\infty(0,T;\dH_\divv^0 (Q_L))$ implies, in turn, 
$u\in C_\omega([0,T];\dH_\divv^0 (Q_L))$. For details, see for instance Lemmas 2.2.3 
and 2.2.5 of \cite{Pok}. 

It holds

\proclaim Lemma 2. 
For any $T>0$ there is a weak solution to $(NS_{f,u_0})$ that satisfies 
$$
{1\over2}{d\over dt}|u(t)|^2_{0,L}+\nu|u(t)|^2_{1,L}\le
\langle f(t),u(t)\rangle_{{-1}, 1; L}.
\eqno(2.7)
$$

\noindent
The proof can be found  for instance in Chapter 3 of \cite{Tem}.  Inequality (2.7) is 
referred to as the {\it (weak) energy inequality} and a weak solution that obeys (2.7) is 
called a {\it Leray-Hopf weak solution}. It is not known if it is unique. If it belongs 
additionally to $L^r(L^s)$ with ${3\over s}+{2\over r}\le1$, $s\in[3,\infty]$ 
(the Ladyzhenskaya-Prodi-Serrin class), it becomes unique and regular, see Serrin \cite{Ser}, Galdi \cite{Gal}, Escauriaza, Seregin, \v Sver\'ak \cite{ES\v S}.

\subsection{2.5. Imbeddings and interpolations}

Firstly, let us show a result needed in Section 4 to estimate the nonlinear term. Let us define
$$
K_{(2.8)}(\alpha,L)
= {2\pi\over L} (2\pi)^{-3}  C_S(1-\alpha) C_S(1),
$$
$$
K_{(2.9)}(\alpha, L)=  \bigg({2\pi \over L}\bigg)^{{ \alpha-2}} C_S\bigg(\alpha-{1\over2}\bigg),
$$
where $C_S (\beta)$ denotes a numerical constant of the optimal $2 \pi$-normalized $3$D Sobolev-Poincar\'e inequality, see subsection 1.2.
\Proposition 1. Assume that $ \alpha\in[0,1] $. Then
$$
| \langle B(a,b), A^\alpha w \rangle_{0,L} | \le K_{(2.8)}(\alpha,L)|a|_{1,L}
|\nabla b|_{L^{3\over2-\alpha}(Q_L)}|A^{\alpha+1\over2} w|_{0,L},
\eqno(2.8)
$$
$$
|v|_{L^{{3\over2-\alpha}}(Q_L)}\le K_{(2.9)}(\alpha, L)
|v|_{\alpha-1,L}^{1\over2}|v|_{\alpha,L}^{1\over2},
\eqno(2.9)
$$
provided the r.h.s.'s are meaningful.

\Proof 
First we perform the estimates for $L={2\pi}$ (where we drop the dependence 
on $Q_{2\pi}$) and next we rescale.

\noindent
{\bf Step 1.} (case $L = {2\pi}$)
The H\"older inequality gives  for $\alpha\in[0,2]$
$$
| \langle B(a,b), A^\alpha w \rangle_{0, 2 \pi}  | = \bigg| (2\pi)^{-d} \int_{-\pi}^{\pi} B(a,b) A^\alpha w   \bigg|  \le (2\pi)^{-d} |a|_{L^6}
|\nabla b|_{L^{3\over2-\alpha} }|A^\alpha w|_{L^{6\over1+2\alpha}}$$
For $\alpha\in[0,1]$ we have
$$
|A^\alpha w |_{L^{6\over1+2\alpha}} \le C_S(1-\alpha)  |A^\alpha w|_{{1-\alpha}, 2 \pi} =  C_S(1-\alpha) |A^{\alpha+1\over2} w|_{0,2\pi},
\eqno(2.10)
$$
where the equality in (2.10) follows from (2.2). Combine the above two estimates to get via the Sobolev-Poincar\'e inequality
$$
| \langle B(a,b), A^\alpha w \rangle_{0, 2 \pi} | \le (2\pi)^{-d}  C_S(1-\alpha) C_S(1) |a|_{1,2\pi}
|\nabla b|_{L^{3\over2-\alpha} }  |A^{\alpha+1\over2} w|_{0,2\pi} \quad \eqno(2.11)
$$
for $ \alpha\in[0,1]$.
Estimate (2.11) is the $Q_{2\pi}$-case of (2.8). Similarly we get $Q_{2\pi}$-case of 
(2.9), namely writing
$$
|v|_{L^{3\over2-\alpha}}\le C_S\bigg(\alpha-{1\over2}\bigg) |v|_{{\alpha-{1\over2}, 2 \pi}} \le C_S\bigg(\alpha-{1\over2}\bigg)
|v|_{\alpha-1,2\pi}^{1\over2}|v|_{\alpha,2\pi}^{1\over2},
\eqno(2.12)
$$
where the later inequality follows from an interpolation, with constant $1$ in view of 
the definition of the norm $|\cdot|_{s,L}$.

\noindent
{\bf Step 2.} (a general $L$ by rescaling) For $h:Q_L\to\R^3$, let us denote its dilation 
$h_{2\pi\over L}\!: Q_{2 \pi}\to\R^3 $ with $\bar h$. Recall from subsection 2.1 that the Fourier coefficients of $h$ and $\bar h$ are identical. Hence 
$|a|_{\beta,L}=|\bar a|_{\beta,2\pi}$, 
$A^\beta w=\big({4\pi^2\over L^2}\big)^\beta
A^\beta \bar w$. One has
$$\eqal{
&| \langle B(a,b), A^\alpha w \rangle_{0, L} | = \bigg({4\pi^2\over L^2}\bigg)^{\alpha + {1 \over 2}} | \langle B(\bar a,\bar b), A^\alpha\bar w \rangle_{0, 2 \pi} | \le \cr 
& (2\pi)^{-d}  \bigg({4\pi^2\over L^2}\bigg)^{\alpha + {1 \over 2}}  C_S(1-\alpha) C_S(1) |\bar a|_{1,2\pi}
|\nabla \bar  b|_{L^{3\over2-\alpha} (Q_{2 \pi})}  |A^{\alpha+1\over2}\bar  w|_{0,2\pi} =: \cr
&  (2\pi)^{-d}  \bigg({4\pi^2\over L^2}\bigg)^{\alpha + {1 \over 2}}  C_S(1-\alpha) C_S(1) \;I,}
$$
where the inequality follows from step 1. In order to scale back $I$ to $Q_L$, we need to know how Lebesgue norms behave under scaling. It holds
$$
|\nabla b|_{L^p(Q_L)}=\bigg({L\over2\pi}\bigg)^{{d\over p}-1}
|\nabla\bar b|_{L^p(Q_{2\pi})}.
$$
Taking this into consideration, we get
$$ I = \bigg({4\pi^2\over L^2}\bigg)^{ {-\alpha }} | a|_{1, L}
|\nabla  b|_{L^{3\over2-\alpha} (Q_{L})}  |A^{\alpha+1\over2}   \omega|_{0, L}. $$
Altogether, the formulas that involve $I$ yield
$$
\big| \langle B(a,b), A^\alpha w \rangle_{0, L}  \big| \le {2\pi\over L} (2\pi)^{-d}  C_S(1-\alpha) C_S(1)  | a|_{1, L}
|\nabla  b|_{L^{3\over2-\alpha} (Q_{L})}  |A^{\alpha+1\over2}   \omega|_{0, L},
$$
which is (2.8). Analogously we get (2.9), because
$$
|v|_{L^{3\over2-\alpha} (Q_L)} = \bigg({L\over2\pi}\bigg)^{{2- \alpha}} |\bar v|_{L^{3\over2-\alpha} (Q_{2 \pi})} \le  \bigg({L\over2\pi}\bigg)^{{2- \alpha}} C_S\bigg(\alpha-{1\over2}\bigg)
|\bar v|_{\alpha-1,2\pi}^{1\over2}|\bar v|_{\alpha,2\pi}^{1\over2}.
$$
\kwadrat

Next, let us present a result  that facilitates the desired scaling-invariance of constants in Corollary 1. To formulate it, we need
\proclaim Definition 2.  $C_I (\alpha, T)$ is a constant of the following $2 \pi$-normalized interpolation inequality
$$
|\nabla f|_{L^4 (0, T; L^{{3\over2-\alpha}}(Q_{2 \pi}))}\le C_I (\alpha, T) |f|^{1 \over 2}_{L^\infty(0,T; |\cdot|_{\alpha, {2 \pi}})} |f|^{1 \over 2}_{L^2(0,T;  |\cdot|_{1+ \alpha, 2 \pi})}.
$$

\noindent
The above interpolation holds for $$f\in L^\infty(0,T; \dH^\alpha (Q_{ 2 \pi}))\cap L^2(0, T; \dH^{1+\alpha}(Q_{ 2 \pi}))$$ in view of
\vskip2pt

\item{(i)} The standard interpolation inequality for integro-differential norms.
\item{(ii)} The Poincar\'e inequality that allows us to write homogeneous Sobolev integro-differential seminorms in the r.h.s. of the interpolation inequality from (i). (Elements of $\dH^\beta$ can be identified with these of  $H^\beta$ that have null mean value).
\item{(iii)} The equivalence of integro-differential and Fourier-based norms.
\vskip2pt
\noindent
Rescaling the above interpolation formula, we obtain
\proclaim Proposition 2.
Assume that  $u\in L^\infty(0,T; \dH^\alpha (Q_L))\cap L^2(0,T; \dH^{1+\alpha} (Q_L))$. Then 
$$
|\nabla u|_{L^4 (0, T; L^{{3\over2-\alpha}}(Q_{L}))}\le C_I (\alpha, T) \bigg({4\pi^2\over L^2}\bigg)^{\alpha-1 \over2} |u|^{{1} \over {2}}_{L^\infty(0,T; |\cdot|_{\alpha, {L}})} |u|^{{1} \over {2}}_{L^2(0,T;  |\cdot|_{1+ \alpha, L})}.
$$

\section{3. Stability}
We are ready to prove Theorem 1. Recall that we work with a given $T_*>0$ and  $\alpha$-strong solution $u$ to $(NS_{f,u_0})$ that exists on $[0, T_*)$ as well as a Leray-Hopf weak solution $v$ to $(NS_{g,v_0})$. The system for the 
difference $w:=u-v$ reads
$$\left\{\eqal{
&w,_{t}+\nu Aw+B(w,u)+B(u,w)-B(w,w)=h\cr
&\divv w=0\cr
&w(0)=u_0-v_0\cr}\right.
\eqno(3.1)
$$
with $h = f-g$.
In subsection 3.1 we will derive higher--order estimates for (3.1) (more precisely, 
$\alpha$-order estimates). Next, we conclude the proof via a blowup argument in subsection 3.2.

\subsection{3.1. Energy estimates}
We are going to test (3.1) with $A^\alpha w$. 

\subsubsection{3.1.1. Admissibility of testing with $A^\alpha w$}
Let us first comment on rigorousness of our estimates. 
We restrict ourselves to the time interval $[0,T_*(g, v_0)\wedge T_*)$,  
where $T_*(g, v_0)$ denotes the blowup time of the $\alpha$-strong solution to $(NS_{g,v_0})$, given by Theorem 2. This solution coincides on the interval $[0, T_*(g, v_0))$ with the interesting for us Leray-Hopf weak solution $v$, again thanks to Theorem 2 (its uniqueness part). Hence $$ w= u-v \in L^2(0, T_*(g, v_0)\wedge T_*; \dH^{1+\alpha}(Q_L)).$$ Consequently
\item{(i)} $A^\alpha w$  is admissible as a test function to  $\nu A w$ in (3.1).
\vskip6pt
\noindent
Next, for $\alpha\ge{1 / 2}$ it holds $$A^{{\alpha \over 2}} w,_{t} \in L^2(0, T_*(g, v_0)\wedge T_*; \dH_\divv^{-1}(Q_L)),$$ thanks to an analogous argument, as the one for formulas (4.9) and (4.10), used for $w$. This and the already known
$$A^{{\alpha \over 2}} w \in L^2(0, T_*(g, v_0)\wedge T_*; \dH^{1}(Q_L)) $$ allows us to write
$$
\langle A^{{\alpha \over 2}} w,_{t}(t), A^{{\alpha \over 2}} w (t)
\rangle_{-1, 1; L}={1\over2}{d\over dt}
|A^{{\alpha \over 2}} w|_{0,L}^2(t).
\eqno(3.2)
$$
Identity (3.2) is the Fourier-series version of the known integro-differential formula for a generalized differentiation of a product, compare for instance Lemma 2.2.5 of \cite{Pok}. Hence
\vskip6pt
\item{(ii)} 
By (3.2) we have justified the admissibility of  $A^\alpha w$ as a test function to the evolutionary part of  (3.1). \vskip6pt
\noindent
Observe that the above justification works well only for 
$\alpha$-order  estimates for $\alpha\ge{1 / 2}$.
Otherwise we do not have sufficient regularity information on $w,_{t}$ to use the duality 
formula (3.2). For instance for $\alpha=0$, one has $u,v\in L^2(\dH^1)$ and 
$u,_{t},v,_{t}\in L^{4\over3}(\dH^{-1})$ (see Lemma 1). In order to have the duality formula (3.2) with 
such low regularity of the time derivative, we would need to assume 
$u,v\in L^4(\dH^1)$ (which is already well within the Ladyzhenskaya-Prodi-Serrin class). 
In this case one can justify the estimates differently, see \cite{CDGG}, proof 
of Theorem 3.3,  in particular pages 61-63. 

\noindent
Finally, 
\vskip6pt
\item{(iii)} 
Testing  the nonlinear and force  terms of (3.1) with  $A^\alpha w$ is admissible. This can be seen directly in the estimates (3.4)--(3.7) 
below.
\subsubsection{3.1.2. Estimates} Testing (3.1) with $A^\alpha w$, we get at a.a. $t\in[0,T_*(g, v_0)\wedge T_*)$
$$\eqal{
&{1\over2}{d\over dt}|A^{\alpha\over2}w|_{0, L}^2+\nu
|A^{\alpha+1\over2}w|_{0,L}^2\le
\cr &| \langle h, A^\alpha w \rangle_{-1,1;L}| 
+ | \langle B(w,u)+B(u,w)-B(w,w), A^\alpha w \rangle_{-1,1;L}|.}
\eqno(3.3)
$$

Let us estimate the force term as follows
$$
| \langle h, A^\alpha w \rangle_{-1,1;L}| = | \langle A^{\alpha-1\over2}h, 
A^{\alpha+1\over2} w  \rangle_{0,0;L}|
\le\varepsilon_2|A^{\alpha+1\over2}w|_{0,L}^2+
{1 \over 4 \varepsilon_2} |A^{\alpha-1\over2}h|_{0,L}^2.
\eqno(3.4)
$$

To control the  nonlinear terms we use (2.8) of Proposition 1 and get
$$\eqal{
I&:=  | \langle B(w,u), A^\alpha w \rangle_{-1,1;L}| \le
K_{(2.8)}(\alpha,L) |w|_{1,L}|\nabla u|_{L^{3\over2-\alpha}(Q_L)}
|A^{\alpha+1\over2}w|_{0,L}\cr
&\le K_{(2.8)}(\alpha,L)|w|_{1-\alpha,L}^{1\over2}
|w|_{1+\alpha,L}^{1\over2}|\nabla u|_{L^{3\over2-\alpha}(Q_L)}
|A^{\alpha+1\over2}w|_{0,L}\cr
&\le K_{(2.8)}(\alpha,L)\bigg({4\pi^2\over L^2}\bigg)^{-{1+\alpha\over4}}
|w|_{1-\alpha,L}^{1\over2}|\nabla u|_{{3\over2-\alpha}(Q_L)}
|A^{\alpha+1\over2}w|_{0,L}^{3\over2}.}
\eqno(3.5)
$$
In (3.5) we use also 
interpolation of $|\cdot|_1$ between $|\cdot|_{1-\alpha}|\cdot|_{1+\alpha}$ 
with constant $1$ (which follows from the definition of $|\cdot |_{s, L} $ and the Cauchy-Schwarz 
inequality) and (2.2). Observe that the term containing $u$ is finite for a.e. $t$ thanks to Proposition 2. Similarly
$$\eqal{
II& :=  | \langle B(u,w), A^\alpha w \rangle_{-1,1;L}|  \cr 
&\le K_{(2.8)}(\alpha,L)K_{(2.9)}
(\alpha,L)|u|_{1,L}|\nabla w|_{\alpha-1,L}^{1\over2}
|\nabla w|_{\alpha,L}^{1\over2}
|A^{\alpha+1\over2}\omega|_{0,L}\cr
&=K_{(2.8)}(\alpha,L)K_{(2.9)}(\alpha,L) \bigg({4\pi^2\over L^2} \bigg)^{1 \over 2} |u|_{1,L}
|w|_{\alpha,L}^{1\over2}| w|_{\alpha+1,L}^{1\over2}
|A^{\alpha+1\over2} w|_{0,L}\cr
&\le K_{(2.8)}(\alpha,L)K_{(2.9)}(\alpha,L)
\bigg( {4\pi^2\over L^2} \bigg)^{1 \over 2} |u|_{1,L}
|w|_{\alpha,L}^{1\over2}\bigg({4\pi^2\over L^2}\bigg)^{-{1+\alpha\over4}}
|A^{\alpha+1\over2} w|_{0,L}^{3\over2}.}
\eqno(3.6)
$$
for the equality above we use '$A^{1 \over 2} = |\nabla|$'  and (2.2) and for the last inequality again (2.2). We begin the estimate of 
the last nonlinear part of (3.3) by invoking (3.6) with $u:=w$
$$\eqal{
III&:=   | \langle B(w,w), A^\alpha w \rangle_{-1,1;L}|  \cr
&\le K_{(2.8)}(\alpha,L)K_{(2.9)}(\alpha,L)
\bigg({4\pi^2\over L^2}\bigg)^{1-\alpha\over4}
|w|_{1,L}|w|_{\alpha,L}^{1\over2}
|A^{\alpha+1\over2}w|_{0,L}^{3\over2}\cr
&\le K_{(2.8)}(\alpha,L)K_{(2.9)}(\alpha,L)
\bigg({4\pi^2\over L^2}\bigg)^{-{\alpha\over2}}
|w|_{1-\alpha,L}^{1\over2}|w|_{\alpha,L}^{1\over2}
|A^{\alpha+1\over2}w|_{0,L}^2,}
\eqno(3.7)
$$
where for the second inequality we interpolate
$|w|_{1,L}\le|w|_{1-\alpha,L}^{1\over2}|w|_{1+\alpha,L}^{1\over2}$ 
and use (2.2). 

Estimates (3.4)--(3.7) plugged into (3.3) yield
$$\eqal{
&{1\over2}{d\over dt}|A^{\alpha\over2}w|_{0, L}^2+\nu
|A^{\alpha+1\over2}w|_{0,L}^2\le
\varepsilon_2 |A^{\alpha+1\over2}w|_{0,L}^2+
{1 \over 4\varepsilon_2} |A^{\alpha-1\over2}h|_{0,L}^2 \cr
&+K_{(2.8)}(\alpha,L)|A^{\alpha+1\over2}w|_{0,L}^{3\over2}\bigg[\bigg({4\pi^2\over L^2}\bigg)^{-{1+\alpha\over4}}
|w|_{1-\alpha,L}^{1\over2}|\nabla u|_{L^{3\over2-\alpha}(Q_L)}
\cr
&+K_{(2.9)}(\alpha,L)\bigg({4\pi^2\over L^2}\bigg)^{1-\alpha\over4}
|u|_{1,L}|w|_{\alpha,L}^{1\over2}\cr
&+K_{(2.9)}(\alpha,L)\bigg({4\pi^2\over L^2}\bigg)^{-{\alpha\over2}}
|w|_{1-\alpha,L}^{1\over2}|w|_{\alpha,L}^{1\over2}
|A^{\alpha+1\over2}w|_{0,L}^{1\over2}\bigg]}
\eqno(3.8)
$$
In (3.8) we need the restriction $\alpha\in\big[{1/2},1\big]$, because we have used 
Proposition 1. Observe that the last summand of (3.8) gives the critically 
growing term
$|A^{\alpha+1\over2}w|_{0,L}^2$. Therefore it may seem more natural to stop 
estimate (3.7) for $III$ at the first inequality and have in consequence the subcritical 
$|A^{\alpha+1\over2}w|_{0,L}^{3\over2}$ instead. Then, however, one needs to deal 
with higher powers of the lower-order-terms. It is possible in case 
$\alpha=1$, but we prefer to keep the energy estimate in the form (3.8) and 
argue for the entire range $\alpha\in\big[{1/2},1\big]$ at once.



In the last-but-one term on the r.h.s. of (3.8) let us use 
$$
|u|_{1,L}=L^{- {d \over 2}} \bigg( {4\pi^2\over L^2} \bigg)^{-{1 \over 2}} |\nabla u|_{L^2(Q_L)}\le
L^{- {d \over 2}} \bigg( {4\pi^2\over L^2} \bigg)^{-{1 \over 2}} L^{ {d (2 \alpha -1) \over 6}}
|\nabla u|_{L^{3\over2-\alpha}(Q_L)},
$$
which follows from (1.1) and the H\"older inequality.  This and  $|f|_{1-\alpha,L}\le|f|_{\alpha,L}$, valid for $\alpha\ge{1/2}$, yields from (3.8) via (2.2)
$$\eqal{
&{1\over2}{d\over dt}|A^{\alpha\over2}w|_{0, L}^2+\nu
|A^{\alpha+1\over2}w|_{0,L}^2\le
\varepsilon_2 |A^{\alpha+1\over2}w|_{0,L}^2+
{1 \over 4\varepsilon_2} |A^{\alpha-1\over2}h|_{0,L}^2 + \cr
&K_{(2.8)}(\alpha,L) 
|A^{\alpha+1\over2}w|_{0,L}^{3\over2} |w|_{\alpha,L}^{1\over2}  |\nabla u|_{L^{3\over2-\alpha}(Q_L)}  \bigg({4\pi^2\over L^2}\bigg)^{-{1+\alpha\over4}}  \bigg[1
+K_{(2.9)} (\alpha,L) L^{ {d (\alpha - 2) \over 3}}  \bigg]\cr
&+K_{(2.8)}(\alpha,L) K_{(2.9)}(\alpha,L)\bigg({4\pi^2\over L^2}\bigg)^{-{\alpha\over2}}|w|_{\alpha,L}
|A^{\alpha+1\over2}w|_{0,L}^2}
$$
Expressing above all the norms of fractional derivatives by the norms of a respective homogenous Sobolev space via (2.2), we arrive at
$$\eqal{
&{1\over2}{d\over dt}|w|_{\alpha, L}^2+\nu {4\pi^2\over L^2}
|w|_{\alpha+1,L}^2\le
\varepsilon_2  {4\pi^2\over L^2}  |w|_{\alpha+1,L}^2+
{1 \over 4\varepsilon_2} \bigg({4\pi^2\over L^2}\bigg)^{-1} |h|_{\alpha-1,L}^2\cr
&+\bigg({4\pi^2\over L^2}\bigg)^{ {3 - \alpha \over 4}} |w|_{\alpha+1,L}^{3\over2}  |w|_{\alpha,L}^{1\over2}  |\nabla u|_{L^{3\over2-\alpha}(Q_L)}   \cr
& \times \bigg({4\pi^2\over L^2}\bigg)^{-{1+\alpha\over4}}  K_{(2.8)}(\alpha,L)   \bigg[1
+K_{(2.9)}(\alpha,L) 
L^{ {d (\alpha - 2) \over 3}}   
\bigg]\cr
&+K_{(2.8)}(\alpha,L) K_{(2.9)}(\alpha,L)\bigg({4\pi^2\over L^2}\bigg)^{1-{\alpha\over2}}|w|_{\alpha,L}
|w|_{1+\alpha,L}^2.}
$$
Let us define
$$\eqal{
&X:={1\over2} |w|_{\alpha,L}^2,\quad &Y:= {4\pi^2\over L^2}
|w|_{\alpha+1,L}^2,\cr
&U:=|\nabla u|_{L^{3\over2-\alpha}(Q_L)}^4,\quad
&H:= |h|_{\alpha-1,L}^2.}
$$
The Young inequality
$$
Y^{3\over4}c(XU)^{1\over4}\le\varepsilon_1 Y+\varepsilon_1^{-3}{27\over256}c^4 XU
$$
used in the third term of the preceding inequality allows us to write
$$
\dX+ (\nu - \varepsilon_1 - \varepsilon_2) Y\le K_4H +K_3XU + K_2X^{1\over2}Y,
\eqno(3.9)
$$
where
$$\eqal{
&K_2= \sqrt{2}K_{(2.8)}(\alpha,L)K_{(2.9)}(\alpha,L)\bigg({4\pi^2\over L^2}\bigg)^{-{\alpha\over2}},\cr
&K_3=\varepsilon_1^{-3}{27\over128} K^4_{(2.8)}(\alpha,L)   \bigg({4\pi^2\over L^2}\bigg)^{-2\alpha-1}  \bigg[1
+K_{(2.9)}(\alpha,L) L^{ {d ( \alpha -2) \over 3}}  \bigg]^4,\cr
&K_4={1\over4 \varepsilon_2} \bigg({4\pi^2\over L^2}\bigg)^{-1}.}
$$
The above choices agree with the definition of $K_2, K_3, K_4$ in subsection 1.2. To see this, consider the formulas for $K_{(2.8)}(\alpha,L), \, K_{(2.9)}(\alpha,L)$ as in subsection 2.5.
The ODI (3.9) will give us stability via a blowup argument.

\subsection{3.2. Proof of Theorem 1 via the blowup argument}

Recall that assumptions of Theorem 1 fix a positive $T$ that satisfies  $T <T_*$, where $T_*$ is the given time of existence of the reference $\alpha$-strong solution $u$.
The proximity assumption (A1) reads
$$\bigg(|u_0-v_0|_{\alpha,L}^2+K_4\int_0^T
|f-g|_{\alpha-1,L}^2 (t) \, dt\bigg) {\rm e}^{ K_3 \int_0^T|\nabla u (t) |_{L^{3\over2-\alpha} (Q_L)}^4 \, dt}<
\bigg({\bar\nu\over K_2}\bigg)^2,
$$
where  $\bar\nu$ is any positive number that satisfies $\bar\nu <\nu -  \varepsilon_1-\varepsilon_2$.

\noindent
{\bf Step 1.} (a lower bound for $T_*(g, v_0)$.)
Let us show that $$T_*(g, v_0)>T.$$ Assume the contrary: $T_*(g, v_0)\le T \;(\;<T_*(f, u_0))$. 
Let
$$
\mu:= \nu -  \varepsilon_1-\varepsilon_2 -\bar\nu,
$$
positive by our assumptions.
Hence the proximity assumption (A1) gives 
$$X(0)<\bigg({\bar\nu\over K_2}\bigg)^2= \bigg({ \nu -  \varepsilon_1-\varepsilon_2 - \mu \over K_2}\bigg)^2.$$ We face now the following  alternative
\vskip6pt

\item{(i)} either $X(t)\le \big({ \nu -  \varepsilon_1-\varepsilon_2 - \mu \over K_2}\big)^2 $ for $t\in[0,T_*(g, v_0))$
\item{(ii)} or $X(t)$ exceeds $\big({ \nu -  \varepsilon_1-\varepsilon_2 - \mu \over K_2}\big)^2$ on 
$[0,T_* (g, v_0))$. Thanks to continuity of $X$ on $[0,T_*(v))$ and the fact that it 
starts below $\big({ \nu -  \varepsilon_1-\varepsilon_2 - \mu \over K_2}\big)^2$, there exists the minimal positive time
$\bar t\in(0,T_*(v))$ such that $X(\bar t)= \big({ \nu -  \varepsilon_1-\varepsilon_2 - \mu \over K_2}\big)^2$.
\vskip6pt

\noindent
Keeping this in mind, observe that ODI (3.9) reads
$$
\dot X+(  \nu -  \varepsilon_1-\varepsilon_2  -K_2X^{1\over2})Y\le K_3XU+K_4H.
$$
It implies for almost any $t \le T_*(g, v_0)$ 
(case (i)) or for almost any  $t\le\bar t$ (case (ii)) that 
$\dot X+\mu  Y\le K_3XU+K_4H$.  Consequently
$$\eqal{
X(t)+\mu\int_0^tY&\le\bigg(X(0)+K_4 \int_0^t H\bigg){\rm e}^{ K_3\int_0^t U(s)ds } \cr
&\le\bigg(X(0)+ K_4 \int_0^TH\bigg) {\rm e}^{ K_3\int_0^T U(s)ds } <\bigg({\tilde\nu\over K_2}\bigg)^2,}
\eqno(3.10)
$$
where the third inequality follows from our assumption (A1).

 In the case (i), we drop the first summand of the l.h.s. of (3.10), so
$\mu\intop_0^tY<\big({\bar \nu\over K_2}\big)^2$ for any $t<T_*(g,v_0)$, hence 
$$\int_0^{T_*(g, v_0)}Y\le {1\over\mu} \bigg({\bar \nu\over K_2}\bigg)^2<+\infty.$$
Since in the case (i) one assumes also that the continuous $X(t)\le\big({\bar \nu\over K_2}\big)^2$ on $[0,T_*(v))$, $T_*(v)$ 
can not be a blowup time.

In the case (ii) we have $X(\bar t)<\big({\bar \nu\over K_2}\big)^2$ from 
(3.10), but $X(\bar t)=\big({\bar \nu \over K_2}\big)^2$ here, which is a contradiction.

As neither (i) nor (ii) can hold, we have contradicted $T_*(g, v_0)\le T$.

\noindent
{\bf Step 2.} 
(proximity estimate) We already know that  $T_*(g, v_0)>T$. Therefore we rewrite the alternative from 
the previous step, plugging there $T$ in place of $T_*(g, v_0)$. Case (ii) is again a contradiction, so (3.10) 
holds, for any $t<T$. We know that $T<T^*(g, v_0)$, so we can let $t\to T$ in (3.10). This 
gives (P1).
\kwadrat

\subsection{3.2. Proof of Corollary 1}
It follows from Proposition 2 in subsection 2.5 and (A1).
\section{4. Regularity}
Here we prove our theorems  on existence of strong solutions to $(NS_{f, u_0})$. 
\subsection{4.1. Proof of Theorem 2}
This theorem serves as an auxiliary result for our main Theorem 1, therefore its proof has been postponed until now. Nevertheless, the approach for proving both Theorem 1 and Theorem 2 is similar. In the former one we had the blowup argument basing on the reference solution $u$. Here, the approximate solution $u^m$ related to $(NS_{f, u_0})$ (defined below) will play the role of a regular reference.
Recall subsection 2.2, for $k_0 \in \dZ$ we denote by 
$$P_{k}:\dH_\divv^0(Q)\to\span\{e^{2\pi ik\cdot{x\over L}}:k\le k_0\}$$
the projection on the low-frequency space.
The approximate solution $u^m$ related to $(NS_{f, u_0})$ is the solution of the ODE
$$\left\{ \eqal{
&u,_{t}^{m}-\nu\Delta u^{m} + P_m(B(u^{m},u^{m})) = P_m f, \cr
&\nabla\cdot u^{m}=0,\cr
&u^{m}(0)=P_m u_0 .\cr} \right.
$$
compare page 57 of  \cite{CDGG}.
As we are already familiar with the proof method and we follow closely the proof of Theorem~3.5 in \cite{CDGG} (except for the step 5),  we omit some details in the considerations below.

\noindent
{\bf Step 1.} (splitting the initial data) 
Let us consider the projection  $P_{k_0}$, whose $k_0$ will be fixed later and decompose the  initial datum $u_0$ into $u_0^{Lo}=P_{k_0}u_0$ and 
$u_0^{Hi}=u_0-P_{k_0}u_0$, whereas the former evolve with the homogenous Stokes, which degenerates in our setting to the homogenous heat system,    i.e.  to
$$\left\{ \eqal{
&u,_{t}^{Lo}+\, \nu Au^{Lo}=0,\cr
&u^{Lo}(0)=u_0^{Lo}.} \right. 
\eqno(4.1)
$$
Observe that for a low frequency data $u_0^{Lo}$ the above Stokes problem admits an 
exact finite Fourier series solution $u^{Lo}$ that belongs to 
$P_{k_0}(\dH_\divv^0)$. Let us choose any $m\ge k_0$ and consider the approximate solution $u^m$ related to $(NS_{f, u_0})$. Then $u^m-u^{Lo}:=u^{m,{Hi}}$ solves
$$\left\{ \eqal{
&u,_{t}^{m,Hi}-\nu\Delta u^{m,Hi} +\cr 
&P_m[(B(u^{m,Hi},u^{m,Hi}))+B(u^{m,Hi},u^{Lo})+B(u^{Lo},u^{m,Hi})]=F,\cr
&\nabla\cdot u^{m,Hi}=0,\cr
&u^{m,Hi}(0)=P_m(u_0^{Hi})=P_m u_0 - P_{k_0} u_0 .\cr} \right.
\eqno(4.2)
$$
with $F=P_mf - P_m[B(u^{Lo},u^{Lo})]$. 

\noindent
{\bf Step 2.} (derivation of an ODI) 
System (4.2) is formally equivalent to (3.1) 
with $P_m B$ in place of $B$ and $w:=u^{m,Hi}$, $u:=u^{Lo}$, $f:=F$. Our `eigenvalue definition' of $A^\alpha$ 
 reduces testing (4.2) with $A^\alpha u^{m,Hi}$ to multiplying  a system 
of ODEs with
$$
\sum_{k\in\dZ,|k|\le m}\lambda_k^\alpha u^{m,Hi}_k \cdot \omega_{L,k},
$$
hence the estimates of subsection 3.1 are justified also for (4.2).
Consequently, along lines of Section 3 we obtain an analogue of the ODI (3.9)
$$
\dX_m+ (\nu - \varepsilon_1 - \varepsilon_2)  Y_m\le K_2X_m^{1\over2}Y_m+K_3X_mU+K_4H_m,
$$
with
$$\eqal{
&X_m:={1\over2} |u^{m,Hi}|_{\alpha,L}^2,\quad &Y_m:= {4\pi^2\over L^2}
|u^{m,Hi}|_{\alpha+1,L}^2,\cr
&U^{Lo}:=|\nabla u^{Lo}|_{L^{3\over2-\alpha}(Q_L)}^4,\quad
&H_m:= |P_m f|_{\alpha-1,L}^2 +K_5 (u^{Lo}) X^{1 \over 2}_m,}
$$
where $$K_5 (u^{Lo}) =  \sqrt{2} \bigg({4\pi^2\over L^2}\bigg)^{\alpha + {1 \over 2} } | 
u^{Lo} \otimes u^{Lo} |_{1+ \alpha, L}.$$
The $K_5 (u^{Lo})$  term of the forcing $H_m$ follows from $P_m [B(u^{Lo},u^{Lo})]$ part 
of $F$. Namely, it holds
$$
|\langle P_m[B(u^{Lo},u^{Lo})], A^\alpha u^{m,Hi} \rangle | \le | {\rm div} (
u^{Lo} \otimes u^{Lo}) |_{\alpha, L} |A^{\alpha}u^{m,Hi} |_{- \alpha, L} $$
Since ${\rm div}$ is the (scalar) multiplication with $k  \big({4\pi^2\over L^2} \big)^{1 \over 2} $ and (2.2) is valid,  we get 
$$
|\langle P_m[B(u^{Lo},u^{Lo})], A^\alpha u^{m,Hi} \rangle | \le \bigg({4\pi^2\over L^2}\bigg)^{{1 \over 2}+\alpha} | 
u^{Lo} \otimes u^{Lo} |_{1+ \alpha, L} |u^{m,Hi} |_{ \alpha, L}, $$
hence the $K_5 (u^{Lo})$ term in $H_m$.
Let us rewrite now our ODI as follows
$$
\dX_m+(\nu - \varepsilon_1 - \varepsilon_2-K_2X_m^{1\over2})Y_m\le X_m [K_3 U^{Lo} + \delta] + {1 \over 4 \delta}K^2_5 (u^{Lo}).
\eqno(4.3)
$$

\noindent
{\bf Step 3.} (gaining a smallness) Since $u^{Lo}$ solves the linear heat system, we control $U^{Lo}$ and $K_5 (u^{Lo})$ in terms of $u^{Lo}_0$.
Thanks to splitting $u^m$ into $u^{m,Hi}$ and $u^{Lo}$ we can now gain smallness of $X_m (0) = {1 \over 2} |P_m u_0 - P_{k_0} u_0|_{\alpha,L}^2$. Namely,
let us fix a positive $\mu < \nu - \varepsilon_1 - \varepsilon_2$ and choose $k_0$ large enough so that for any $m\ge k_0$ holds
$$
\nu - \varepsilon_1 - \varepsilon_2-K_2X_m^{1\over2}(0)\ge\mu\Leftrightarrow X_m(0)\le
\bigg({\nu - \varepsilon_1 - \varepsilon_2-\mu\over K_2}\bigg)^2.
\eqno(4.4)
$$

\noindent
{\bf Step 4.} ($m$-uniform $\alpha$-regularity bound) Fix any $\sigma \in (0,1)$.
By time continuity of $X_m$, there exists $T_m$ such that 
$\nu - \varepsilon_1 - \varepsilon_2-K_2X_m^{1\over2}(t)\ge{\sigma \mu}$ 
for $t\in[0,T_m]$. It implies in (4.3)
$$
\dX_m+ \sigma Y_m\le  X_m [K_3 U^{Lo} + \delta] + {1 \over 4 \delta}K^2_5 (u^{Lo})
$$
for  $t\le T_m$,   i.e. 
$$
X_m(t) + {\sigma} \int_0^t Y_m\le\bigg(X_m(0)+ {1 \over 4 \delta} \int_0^t  K^2_5(u^{Lo}) \, d \tau \bigg)e^{ \int_0^t  (K_3 U (\tau) + \delta) d \tau}.
\eqno(4.5)
$$
Since in view of (4.4) $ X_m(0)\le
\big({\nu - \varepsilon_1 - \varepsilon_2-\mu\over K_2}\big)^2$, there exists $T_0>0$ such that r.h.s. of (4.5) stays below $\big({\nu - \varepsilon_1 - \varepsilon_2-{\sigma \mu} \over K_2}\big)^2$ for $t\le T_0$ small 
enough. Consequently
$$
X_m(t) \le \bigg({\nu - \varepsilon_1 - \varepsilon_2-{\sigma \mu} \over K_2}\bigg)^2 \Leftrightarrow \nu - \varepsilon_1 - \varepsilon_2-K_2X_m^{1\over2}(t)\ge {\sigma \mu}\eqno(4.6)$$
for $t\le T_0$ and independently from $m$.

Using (4.6) in (4.3) allows to conclude that (4.5) holds for $t\le T_0$ uniformly in $m$. Thus we have an 
additional $L^\infty(\dH^\alpha)\cap L^2(\dH^{1+\alpha})$ estimate for 
$u^{m,H}$, hence $u^m$, uniform for $m\ge k_0$.

\noindent
{\bf Step 5.} (time-continuity $C(\dH^\alpha)$). In this step we do not need a precise control over constants. Conseqently $C$ denotes a general constant, that may vary between lines.
Let us divert from \cite{CDGG} and use the following duality estimate for the 
weak solution (2.5) to $(NS)_{f, u_0}$
$$\eqal{
&\big|\int_0^T\langle u,_{t},
\varphi\rangle_{{-1},1;L}\big|\cr
&\le C \int_0^T \big(
|f|_{{\alpha-1,L}}|\varphi|_{{1-\alpha,L}}+
|\nabla u|_{\alpha,L}|\nabla\varphi|_{{-\alpha},L}+
| \divv(u\otimes u)|_{\alpha-1,L}|\varphi|_{1-\alpha,L} \big).}
\eqno(4.7)
$$
The last summand above can be estimated by
$$\eqal{ 
&\int_0^T
| u\otimes u|_{\alpha,L}|\varphi|_{1-\alpha,L} 
\le C \bigg(\int_0^T \int_{Q_L}|A^{\alpha \over 2} u|^{10\over3}\bigg)^{3\over5}
\bigg(\int_0^T \int_{Q_L}|u|^5\bigg)^{2\over5}|\varphi|_{L^2(|\cdot|_{1-\alpha}) }\cr}
\eqno(4.8)
$$
using H\"older inequality.
Plugging (4.8) into (4.7) gives
$$\eqal{
&\big|\int_0^T\langle u,_{t},
\varphi\rangle_{{-1},1;L}\big| \le \cr
&C |\varphi|_{L^2(|\cdot|_{1-\alpha}) } \big[|f|_{L^2(|\cdot|_ {\alpha-1})}+
|u|_{L^2( |\cdot|_{1+\alpha})} +|A^{\alpha \over 2} u|_{L^{10\over3}(L^{10\over3})}
|u|_{L^5(L^5)}^2\big].\cr}
\eqno(4.9)
$$
The norms on the r.h.s. of (4.9)   are finite thanks to our assumption related to $f$ and
to the $L^\infty(\dH^\alpha)\cap L^2(\dH^{1+\alpha})$ estimate from previous steps for $u$. 
Specifically, the last summand in (4.9) is finite by parabolic embedding ($L^{10\over3}$ 
norm) and  by interpolation 
$$L^\infty(\dH^\alpha)\cap L^2(\dH^{1+\alpha})\hookrightarrow 
L^5(\dH^{\alpha+{2\over5}})\hookrightarrow L^5(L^5)$$ for $\alpha\ge{1/2}$. 
Hence (4.9) means that
$$
u,_{t}\in(L^2(\dH_\divv^{1-\alpha}))^*=L^2(\dH_\divv^{\alpha-1}). \eqno(4.10)
$$
This information interpolated (in the sense of `espaces des traces', see for instance 
Lemma 2.2.4 in \cite{Pok}) with 
$u\in L^2(H_\divv^{\alpha+1})$ yields $u\in C(\dH^\alpha)$.

\noindent
{\bf Step 6.} (uniqueness) 
The $L^\infty(\dH^\alpha)$ regularity implies for $\alpha\ge{1/2}$ that 
we are in the Prodi-Serrin class, where Leray-Hopf 
solutions are unique.

\noindent
{\bf Step 7.} (blowup criterion) Assume on the contrary that $T_*(u_0,f)$ is the maximum existence time and at the same time
$$
\max_{t \le  T_*(u_0,f)} |u(t)|_{\alpha,L} <\infty
$$
then, by definition of the $C(I; V)$-Bochner norm, $u \in C([0,  T_*(u_0,f)] ; \dH^\alpha)$. In particular, we can restart the evolution from $u (T_*(u_0,f)) \in \dH^\alpha$ and in view of steps 1--6, there exists $T_1 > T_*(u_0,f)$ and the unique (in Leray-Hopf class)
$$C([ T_*(u_0,f),T_1],\dH_\divv^\alpha)\cap L^2(T_*(u_0,f),T_1,\break\dH_\divv^{1+\alpha})$$
solution to $(NS_{f,u(T_*(u_0,f))})$. It satisfies for 
a.a. $t\in[0,T_1]$ the weak formulation (2.6). Therefore it is a weak solution on $[0,T_1)$ to $(NS_{f,u_0})$. 
Hence $T_*(u_0,f)$ is not the maximal existence time.
\kwadrat
\subsection{4.2. Proof of Lemma 1}
Here we prove the caloric lower bound for $T_*(f, u_0)$. From (4.5) in the step 4 of the proof of Theorem 2, we see that any ${T_0}>0$ that yields
$$\eqal{
&\bigg(\bigg({\nu - \varepsilon_1 - \varepsilon_2-\mu\over K_2}\bigg)^2+ {1 \over 4 \delta} \int_0^{T_0}  K^2_5(u^{Lo} (\tau)) \, d\tau \bigg)e^{ \int_0^{T_0}  (K_3 U (\tau) + \delta) d \tau} \cr
&\le \bigg({\nu - \varepsilon_1 - \varepsilon_2-{\sigma \mu} \over K_2}\bigg)^2 }\eqno(4.11)
$$
provides the $m$-independent bound 
$$
\sup_{t \le {T_0}}  X_m(t) + {\sigma} \int_0^{T_0} Y_m\le \bigg({\nu - \varepsilon_1 - \varepsilon_2-{\sigma \mu} \over K_2}\bigg)^2.  \eqno(4.12)
$$
Let us reformulate (4.11) to obtain
$$\eqal{
& \bigg( {K_2 \over {\nu - \varepsilon_1 - \varepsilon_2-{ \sigma \mu} }} \bigg)^2 {1 \over 4 \delta} \int_0^ {T_0}  K^2_5(u^{Lo}) \, d\tau  \le \cr
& e^{ -\int_0^ {T_0}  (K_3 U (\tau) + \delta) d \tau} - \bigg({{(1-\sigma)\mu} \over{\nu - \varepsilon_1 - \varepsilon_2-{\sigma \mu}} }-1\bigg)^2, }
$$
which clarifies (A3), after one takes into account the formulas 
 $$K_5^2 (u^{Lo}) =  2 \bigg({4\pi^2\over L^2}\bigg)^{2 \alpha +1} | 
u^{Lo} \otimes u^{Lo} |^2_{1+ \alpha, L}, $$
$$U^{Lo} :=|\nabla u^{Lo} (\tau)|_{L^{3\over2-\alpha}(Q_L)}^4,$$
as in  the step 4 of the proof of Theorem 2. The remaining to prove bound follows from (4.12).
\kwadrat

\subsection{4.3. Proof of Theorem 3}

Currently we find ourselves in an easier situation than when proving Theorem 2, because splitting the initial data to gain smallness is unnecessary - a smallness is already assumed. 
Hence we get for the Fourier approximations $u^m$
$$
\dX_m+ (\nu - \varepsilon_1 - \varepsilon_2 - K_2X_m^{1\over2})  Y_m\le K_4H,
$$
with
$$
X_m:={1\over2} |u^{m}|_{\alpha,L}^2,\quad Y_m:= {4\pi^2\over L^2}
|u^{m}|_{\alpha+1,L}^2, \quad H:= | f|_{\alpha-1,L}^2,
$$
compare the computations that provided (4.3) in step 2 of the proof of Theorem 2. We finish our proof via a blowup argument, analogously to the proof of our stability result, compare subsection 3.1.
In fact Theorem 3 can be seen also as a stability 
result with null initial data and null reference solution $u$ in $(A1)$.
\kwadrat

\section{5. Concluding remarks}
Let us recall from our introduction that the most ambitious task related to stability studies of Navier-Stokes is to obtain regularity {\it in the large} by stability. Some of the ideas how to complete this task are related to so called 
  {\it schemes for numerical verification of regularity}. They were first presented in Chernyshenko, Constantin, Robinson \& Titi [CCRT], further generalized in Dashti \& Robinson \cite{DR} and refined in Mar\'\i n-Rubio, Robinson \& Sadowski \cite{M-RRS}. Some remarks on these schemes follow, in relation to our results.
  \subsection{5.1. {A posteriori} numerical verification of regularity}
Let us consider an approximate solution to unforced $ (NS_{0,u_0})$. To fix ideas, let the superscript $\cdot^n$ denote the projection of $u_0$ on $ \omega_{L,k}$, $k = 1, \dots, n $. The  $n$-th approximation starts at $u^n_0$. It is smooth and satisfies 
$$\left\{\eqal{
&u,^n_{t}+u^n\cdot\nabla u^n-\nu\Delta u^n+\nabla p^n=u^n\cdot\nabla u^n - (u^n\cdot\nabla u^n)^n\cr
&\divv u^n=0\cr
&u^n(0)=u^n_0\cr}
\right.
$$
Using stability formula  $(A1)$ for $u_n$ as the reference smooth solution, we have that $u$ is smooth, provided
$$
\bigg(|u^n_0-u_0|_{\alpha, L}^2+ K_4 \int_0^T 
|u^n\cdot\nabla u^n - (u^n\cdot\nabla u^n)^n|^2_{\alpha-1, L} (t) \, dt\bigg) \; {\rm e}^{ K_3 \int_0^T
|\nabla u^n (t)|_{L^{3\over2-\alpha} (Q_L)}^4 \!\! dt} 
<{\bar\nu^2\over K_2^2}
\eqno(C)
$$
Hence to conclude that $u$ is regular, one needs to find $n$ such that condition  (C) holds. The term involving initial condition vanishes, but the rest is troublesome. Up to now, one copes with them in the case $\alpha = {1 \over 2}$ and $u_0 \in H^1$ (observe a mismatch in regularity for data and for stability) {\it assuming a priori that $u$ is regular}, which implies that $$ \int_0^T 
|u^n\cdot\nabla u^n - (u^n\cdot\nabla u^n)^n|^2_{\alpha-1, L} (t)  \, dt $$ vanishes and that $${  \int_0^T
|\nabla u^n (t)|_{L^{3\over2-\alpha} (Q_L)}^4 \, dt} $$ is bounded; see Theorem 6 and Lemma 7 in  \cite{M-RRS}. This need for {\it a priori} assumption of smoothness of $u$ is the main difficulty of this approach.\\
\noindent
  \subsection{5.2. Further research}
Let us finally make a few observations. \\
(i) It may be interesting to try to match the stability condition $(C)$ {\it ad hoc} numerically, where our scaling may be helpful to speed up the computations by reshaping constants (recall formulas for $K_3$ and $K_4$). {\it Ad hoc numerically} means here that one checks whether (C) holds for some $n$ by a numerical scheme. Consequently no assumption on regularity of $u$ is needed, but there is no evidence by now that (C) can be verified numerically. \\
(ii) Even more interesting would be to falsify  numerically the regularity assumption of $u$ as follows: Assume that $u$ starting from $u_0$ is regular and, using this regularity, provide a convergence rate of $$\int_0^T 
|u^n\cdot\nabla u^n - (u^n\cdot\nabla u^n)^n|^2_{\alpha-1, L} (t)  \, dt$$ to zero and the upper bound on $${  \int_0^T
|\nabla u^n (t)|_{L^{3\over2-\alpha} (Q_L)}^4 \, dt}. $$ This implies validity of (C) for  $n \ge n_0$. If the numerics showed otherwise, it would indicate at lack of regularity. \\
(iii) Our introduction of more general $\alpha$ should allow to match the stability condition and the initial datum space with a certain $\alpha < 1$. It is the object of our current studies. \\

\section{References}

\parindent 15mm

\item{[ADT]} 
Auscher, P., Dubois, S., Tchamitchian, P.,
{\it On the stability of global solutions to Navier-Stokes equations in the space}, J. Math. Pures Appl. (9) 83 (2004), no. 6, 673--697. 

\item{[BLNNT]} 
Bardos, C., Lopes Filho, M. C., Niu, D., Nussenzveig Lopes, H. J.; Titi, E. S., {\it Stability of two-dimensional viscous incompressible flows under three-dimensional perturbations and inviscid symmetry breaking}. SIAM J. Math. Anal. 45 (2013), no. 3, 1871--1885.

\item{[BdV]}
Beir\~ao da Veiga, H., Secchi, P., {\it $L^p$-stability for the strong solutions of the Navier-Stokes equations in the whole space}. Arch. Rational Mech. Anal. 98 (1987), no. 1, 65--69.

\item{[CDGG]} Chemin, J.-Y., Desjardins B., Gallagher, I., Grenier, E.,
{\it Mathematical geophysics}, Oxford 2006.

\item{[CCRT]} Chernyshenko, S., Constantin, P., Robinson, J., Titi, E., 
{\it A posteriori regularity of the three-dimensional Navier-Stokes  from 
numerical computations}, Com. Math. Phys. 48 (2007), no. 6, 065204.

\item{[CF]} Constantin P., Foias, C., {\it Navier-Stokes Equations}, Chicago 
1988.

\item{[DR]} Dashti, M., Robinson J., {\it An A Posteriori Condition on the Numerical Approximations of the Navier--Stokes Equations for the Existence of a Strong Solution}, SIAM J. Numer. Anal., 46(6), 3136--3150.

\item{[ES\v S]} Escauriaza, L., Seregin, G., \v Sver\'ak, V., 
{\it $L^{3,\infty}$-solutions of\break Navier--Stokes equations and backward 
uniqueness}, Uspekhi Mat. Nauk 58 (2003), no 2 (350), 3--44.

\item{[Gal]} Galdi, G., {\it An introduction to the Navier-Stokes initial 
boundary value problem}, Birkh\"auser, 2000.

\item{[I]} Iftimie, D., {\it The 3D Navier-Stokes equations seen as a perturbation of the 2D Navier-Stokes
equations}, Bull. Soc. Math. France, 127 (1999), 473--517.

\item{[KT]} Koch, H., Tataru, D., {\it Well-posedness for the Navier-Stokes equations} Adv. Math. 157 (2001), no. 1, 22 -- 35.

\item{[M-RRS]} Mar\'\i n-Rubio, P., Robinson, J., Sadowski, W., {\it Solutions of the 
3D Navier-Stokes equations for initial data in $\dH^{1\over2}$: robustness of 
regularity for bounded sets of initial data in $\dH^1$}, J.M.A.A. 400 (2013), 
no 1, 76--85.

\item{[M1]} Mucha, P.,  {\it  Stability of 2D incompressible flows in $R^3$}, J. Differential Equations 245 (2008), no. 9, 2355--2367.

\item{[M2]} Mucha, P.,  {\it Stability of nontrivial solutions of the Navier-Stokes system on the three dimensional torus}, J. Differential Equations 172 (2001), no. 2, 359--375

\item{[Pok]} Pokorn\'y, M., {\it Navier-Stokes system}. Prague, 2010, available online at \\ \url{http://ssdnm.mimuw.edu.pl/pliki/wyklady/pokorny.pdf}, acc. 21.01.2016.

\item{[Ser]} Serrin, J., {\it The initial value problem for the Navier-Stokes 
system}, Nonlinear Problems (Proc. Sympos., Madison, Wis., 1962), 69 -- 98, Univ. of Wisconsin Press, Madison, Wis., 1963.

\item{[Tem]} Temam, R., {\it Navier-Stokes Equations and Nonlinear Functional 
Analysis}, Second edition, SIAM, 1995.

\item{[Za]} Zaj\c aczkowski, W. M., {\it Some global regular solutions to 
NSE},\break Math. Meth. Appl. Sc. 30 (2006), 123--151.

\item{[ZaZa]} Zaj\c{a}czkowski, W. M.,  Zadrzy\'nska, E., {\it Global 
regular solutions with large swirl to the NSE in a cylinder}, J. Math. Fluid 
Mech. 11 (2009), 126--169.

\item{[Z]} Zeidler, E., {\it Nonlinear Functional Analysis and its Applications II/A}, Springer, 1990.

\bye


$$
K_{(2.7)}(\alpha,L)(\alpha,L)=\cases{c_s(1-\alpha)(1+c_p(1-\alpha)L^{1-\alpha})
\big({4\pi^2\over L^2}\big)^{\alpha-1\over2}& for\ \ $\alpha<1$\cr
1&for\ \ $\alpha=1$\cr}
\eqno(2.8)
$$

$$
|v|_{L^{3\over2-\alpha}}\le\tilde K_{(2.9)}|v|_{\dH^{\alpha-1}}^{1\over2}
|v|_{\dH^\alpha}^{1\over2}\quad {\rm for}\ \ \alpha\in\bigg[{1\over2},2\bigg)
\eqno(2.9)
$$
\bye